\documentclass[9pt, twoside]{article}
\usepackage{japs-ISOSS}
\usepackage{float}
\usepackage{bbm}
\usepackage{amsfonts}
\usepackage{mathrsfs}
\usepackage{amsmath}
\usepackage{dsfont}
\usepackage{amssymb}
\usepackage{graphicx}

\begin{document}
\thispagestyle{empty} \setcounter{page}{1} \setcounter{section}{0}
\setcounter{equation}{0}\setcounter{theorem}{0}\setcounter{footnote}{0}
\setcounter{table}{0} \setcounter{figure}{0}
\numberwithin{equation}{section}





\title{A New Algorithm to Simulate the First Exit Times of a Vector of Brownian
Motions, with an Application to Finance}

\author{Chiu-Yen Kao} 

\address{Department of Mathematical Sciences, Claremont McKenna College, USA. This author is partially supported by NSF DMS-1318364.
}

\email{Ckao@claremontmckenna.edu} 

\author{Qidi Peng}

\address{Institute of Mathematical Sciences, Claremont Graduate University, USA}

\email{Qidi.Peng@cgu.edu}

\author{Henry Schellhorn}

\address{Institute of Mathematical Sciences, Claremont Graduate University, USA}

\email{Henry.Schellhorn@cgu.edu}

\author{Lu Zhu}

\address{Department of Accounting and Finance, University of Wisconsin,  Eau Claire, USA}

\email{zhul@uwec.edu}

\vspace{.2in}

\centerline{\sc{summary}} \vspace{0.1in} \baselineskip.20in
\centerline{
\begin{minipage}{4.5in} \small{ \noindent
We provide a
new methodology to simulate the first exit times of a vector of
Brownian motions from an orthant.  This new approach can be used to simulate the first exit times of dimension higher than two. When at least one Brownian motion has non-zero drift, the joint density function of the first exit times in $N$ dimensions
needs to be known, or approximated. However, when the drifts are all zero, a simpler simulation algorithm is obtained without using the joint density function.
}\end{minipage}}

\vspace{.15in} \baselineskip.18in


\centerline{ \begin{minipage}{4.5in} \small{ \noindent {\it Keywords
and phrases:} First exit times; Correlated Brownian motions; Gaussian copula; Multiple roots transformation}\end{minipage}}

\vspace{.1in}

\centerline{ \begin{minipage}{4.5in} \small{ \noindent {\it AMS
Classification:} 60G15; 60E05; 65C05; 65C30}\end{minipage}}


\section{Introduction}
\label{Introduction}

There are two main approaches to calculate the expected value of a function
of joint first exit times of a process $X$ from a sufficiently regular domain in $%
\mathbb{R}^{N}$, whether bounded or unbounded. The first approach is to set
up the problem as a partial differential equation. When the process is
 Brownian motion, this equation is of either elliptic or parabolic type. The second approach is Monte Carlo
simulation. The simplest implementation of Monte Carlo method is to discretize the stochastic differential
equation of the underlying process, using for instance the Euler scheme, and
to simulate the process until it reaches the barrier, i.e., the boundary of
the domain (we will discuss and apply this method in Section \ref{K-S Test} and Section \ref{Applications}). We refer to this indirect method as Euler-based
Monte Carlo. According to  Fahim et al. \cite{Fahim}, it is well known that numerical methods based on partial
differential equations suffer from the curse of dimensionality: calculations
and memory requirements typically increase exponentially as the dimension $N$ increases. This is not the case for Monte
Carlo simulation, where calculations are in general independent of the dimension of the
problem.

A hybrid class of approaches uses Monte Carlo method to simulate the first exit times
directly, i.e., without first simulating the underlying process $X$. There
are two main advantages of this direct simulation approach, compared to
Euler-based Monte Carlo. The first one relies on the fact that, in many applications, exit
from a bounded domain is a rare event. It may take a very large amount of
time (sometimes an infinite amount of time) for a process $X$ to reach the
boundary. The traditional solution to this problem is to use importance
sampling. For advantages and pitfalls of importance sampling in financial
engineering, we refer to \cite{Glasserman} and the references
therein. Another solution, in case the importance sampling density is
difficult to construct, is to use interacting particle systems as in \cite{Carmona}. However, we can avoid this problem with direct simulation of the first exit times. The second advantage is its accuracy compared to Euler-based Monte Carlo method. Euler-based Monte Carlo tends to overestimate the first exit time: since the process is discretized, the method misses an excursion, i.e., an outcome where the process is in the domain at two contiguous epochs, but outside the domain at some time in-between these epochs. We note however that there are some
techniques by Broadie et al. \cite{Broadie}, Huh and Kolkiewicz \cite{Huh} and
Shevchenko \cite{Shevchenko} to alleviate this problem.

It is worth noting that direct first exit time Monte Carlo algorithms are not plagued by these two problems mentioned above but relying on the following two principles. First, the joint density of
the first exit times must be either known or approximated correctly. This problem
is often solved by a partial differential equation method, and is quite a
difficult analytical problem in high dimensions if the different components
of the process $X$ are correlated. Secondly, exit times must be sampled from
the joint distribution of the corresponding processes without, if possible, resorting to inversion and
conditional Monte Carlo, which are very computationally intensive. For
bounded domains, there is already rich literature. In \cite{Muller} the
bounded domain is approximated by a series of spheres; in \cite{Milstein} the bounded domain is approximated by a series of parallelepipeds. When the joint density of first exit times for general diffusions
are unknown, the latter authors approximate the diffusion $X$ by a
 Brownian motion on small bounded domains. We note that these methods provide better approximation than the simple Euler-based Monte Carlo method does, since the barrier does not need to be corrected.

Simulating the first exit time from an unbounded domain is a slightly
different problem. We focus on a particular application in finance: the
simulation of the joint default times of a set of $N$ obligors. Black and Cox \cite{Black} were the first to model the default time of an
obligor as the first time the value of the asset returns $X$ crosses a
lower, constant barrier. There is however no upper barrier on $X$. There have
been many variations on this first exit time approach (which the finance
literature refers to often as first passage time models). In \cite{Black}, $X$ is  a Brownian motion. Collin-Dufresne and Goldstein \cite{Collin-Dufresne}
use a more sophisticated model with mean-reversion in credit spreads. In the last ten years
the interest has shifted, from examining a single obligor, to examining a
large collection of obligors. Pools of large obligors were the building
blocks of the credit default swaps which were heavily
traded before the subprime  crisis of $2008$. It is now recognized that part of
the blame of this crisis should be ascribed to the poor risk management of
these contracts. Indeed, the task of estimating the parameters of these
high-dimensional processes and of simulating them is still recognized as a
difficult problem (see \cite{Anonymous}). The difficulty of the simulation part of the
problem does not come from the complexity of the unbounded domain (usually,
a plain orthant), but the complexity of simulating in high dimensions. Zhou \cite{Zhou} provides a formula for the expected first default time (i.e., the
expected value of the minimum of the first exit times), which is valid only
when $N=2$. Metzler \cite{Metzler2} derives the formula for the joint
density of the first exit times when $N=2$, based on earlier work by Iyengar \cite{Iyengar}. We generalize this formula to take into account
the cases that were not given in \cite{Metzler2}. Metzler \cite{Metzler1,Metzler2} also develops an algorithm to
simulate the first exit times when the underlying two-dimensional Brownian motion has zero drift. We remind the readers that this is a more general
problem than merely finding the expected value of the first default time, as
the price of several complicated derivatives cannot be described as a
function of the latter. Compared to Metzler \cite{Metzler1,Metzler2}, our algorithm can be applied in the following cases:
                        \begin{enumerate}
                          \item $N=2$ and at least one Brownian motion has non-zero drift.
                          \item $N>2$ and all the Brownian motions are with zero drifts.
                        \end{enumerate}
                         The main idea of generating the vector of first exit times in
this paper is to transform the joint first exit times to some joint Chi-squared variables
and then to generate the latter vector starting from standard Gaussian vector. From
 (\ref{transform}) below we see that, the transformations to Chi-squared distribution have single
roots when the drift parameter $\mu_i=0$ and multiple roots when $\mu_i\neq0$. Hence our
main results should be stated under two different scenarios: with drift and
without drift. It is worth noting that, when the Brownian motions are correlated, our algorithm is not exact (except for two-dimensional Brownian motion with zero drift), however, statistical test and
empirical study show that the error of simulation is acceptable in practice (see
Sections \ref{K-S Test} and 4).

The remainder of the paper proceeds as follows. In Section 2, we present our main
theoretical results. The simulation of first exit times becomes the simulation
of some 'basic' random variables. When the drift vector of X is non-zero the
'basic' random variables consist of uniform variables and Gaussian variables,
and the simulation method depends on an explicit representation of the joint
distribution of the first exit times, which is known only when $N=2$. When the
drifts are all zero, the method becomes simpler and much more powerful: the 'basic' random
variables are Gaussian, and the explicit representation of the joint distribution
of the first exit times is not needed. Section 3 summarizes the algorithms. In
Section 4, we present an application of our methodology to the first exit time
problem in finance. We compare our results to Zhou's analytical results \cite{Zhou} when
$N=2$, and to Euler-based Monte Carlo simulation when $N\ge2$.

\section{Main Results}

\subsection{Statistical Modeling}

For some integer $N\geq 2$, let us consider $\{(X_{1}(t),\ldots ,X_{N}(t))\}_{t\ge0}$,
an $N$- dimensional  Brownian motion, i.e. for $k\in \{1,\ldots
,N\}$,
\begin{equation}\label{dXk}
\,\mathrm{d}X_{k}(t)=\mu _{k}\,\mathrm{d}t+\sigma _{k}\,\mathrm{d}W_{k}(t),
\end{equation}%
where $\{W_{k}(t)\}_{t\geq 0}$'s are standard Brownian motions and $\mu
_{k}\in \mathbb{R},\sigma _{k}>0$ are  drift and volatility parameters respectively.

To introduce our estimation method let us start by some notation conventions:

\begin{itemize}
\item Denote by $(X_{1}(0),\ldots ,X_{N}(0))=(x_1,\ldots,x_N)\in\mathbb{R}^N$ the
initial values.
\item Let $(b_1,\ldots,b_N)\in\mathbb{R}^N$ be a vector of barriers. We assume that $b_k\neq x_k$ for all $%
k\in\{1,\ldots,N\}$.

\item Let $(\tau_1,\ldots,\tau_N)$ be the vector of first exit times of $(X_1(t),\ldots,X_N(t))$ to the
barriers $(b_1,\ldots,b_N)$. More precisely, for $k\in\{1,\ldots,N\}$,
$$
\tau_k=\inf\big\{t>0:X_k(t)=b_k,~X_k(0)=x_k\big\}.
$$
\end{itemize}

We are mainly interested in the problem of generating the joint first exit
times $(\tau_1,\ldots,\tau_N)$. Let us recall that its marginal density is well known in the literature: let
$\{X(t)\}_{t\geq0}$ be a  Brownian motion with drift $\mu\in%
\mathbb{R}$ and variance $Var(X(t))=\sigma^2t$ with $\sigma>0$ and let $b\in%
\mathbb{R}$ be some given barrier and $\tau_b$ be the waiting time until
the process $\{X(t)\}_{t\ge0}$ first exits the horizontal line $Y=b$. Conditional on $X(0)=x\neq b$ and $\frac{\mu}{b-x%
}\geq0$, $\tau_b$ has the probability density
function (see for example \cite{Iyengar}):
\begin{equation}  \label{densityT}
f_{\tau_b}(t)=\frac{|b-x|}{\sigma\sqrt{2\pi t^3}}\exp\Big(-\frac{(\mu
t-(b-x))^2}{2\sigma^2 t}\Big).
\end{equation}
The function $f_{\tau_b}$ is in fact an inverse Gaussian density. Remark that, if $\frac{\mu}{b-x}%
<0$, this density of first exit time is defective, i.e., $\mathbb{P}%
(\tau_b<+\infty)<1$. To overcome this inconvenience, we
can generalize the density of first exit time as a mixture of distributions:
\begin{equation*}
g_{\tau_b}(t)=f_{\tau_b}(t)\mathds{1}_{[0,+\infty)}(t)+\Big(%
1-\int_0^{+\infty}f_{\tau_b}(u)\, \mathrm{d} u\Big)\delta_{+\infty}(t),
\end{equation*}
where $\mathds{1}_{[0,+\infty)}$ denotes the characteristic function and $%
\delta_{+\infty}$ is the Dirac measure.

A simple componentwise simulation could be used when the Brownian motions are independent. When the vector of Brownian motions consists of correlated components, the problem of simulating the joint first exit times could be rather challenging and we mainly focus on this problem.

\subsection{First Exit Times of  Brownian Motions with General
Drifts and Multiple Roots Transformation}
\label{Drifted}  Generally, we suppose that $\{X_1(t),\ldots,X_N(t)\}$ are
correlated drifted Brownian motions with $|Corr\big(X_i(t),X_j(t)\big)|<1$, for $i,j\in\{1,\ldots,N\}$ and $i\neq j$. Using the property of inverse Gaussian vector, we introduce a transformation of $
\tau_i^{\prime }$'s to Chi-squared random variables. More precisely, for $i\in\{1,\ldots,N\}$, define the transformation $H_i$ as:
\begin{equation}  \label{transform}
H_i(\tau_i)=\frac{(\mu_i \tau_i-(b_i-x_i))^2}{\sigma_i^2\tau_i}\sim
\chi^2(1),
\end{equation}
where $\chi^2(1)$ denotes the Chi-squared distribution with one degree
of freedom.

Now we are going to state the main results of this part. First we suppose $
\mu$ takes any real value. The following proposition extends (\ref{transform}) to high dimension.
\begin{proposition}
\label{prop:simulation}
 (See \cite{Shuster}) Let $(\tau_1,\ldots,\tau_N)$ be the vector of first exit times
of an $N$-dimensional  Brownian motions with parameters $\big(%
\mu_i,\sigma_i,x_i,b_i\big)$, $i\in\{1,\ldots,N\}$. Then there exists a vector of Chi-squared random variables with one degree of freedom $\big(
\chi_1^2,\ldots,\chi_N^2\big)$ such that the following relation holds:
\begin{equation*}
\big(H_1(\tau_1),\ldots,H_N(\tau_N)\big)\sim \big(
\chi_1^2,\ldots,\chi_N^2\big),
\end{equation*}
where $\{H_i\}_{i=1,\ldots,N}$ are given in (\ref{transform}).
\end{proposition}
Note that through the remaining part of
this paper we suppose that for all $i\in\{1,\ldots,N\}$, $\frac{\mu_i}{
b_i-x_i}\ge 0$, so that the density of $\tau_i$ is not defective. Now we explain our idea of simulation. By Proposition \ref{prop:simulation}, we may consider using an inverse transform to simulate $(\tau_1,\ldots,\tau_N)$, however this inverse transform has multiple roots.  Hence, by
using a similar but more generalized method in \cite{Michael}, we are able to
generate the first exit times approximately starting from the joint
distribution of $(\chi_1^2,\ldots,\chi_N^2)$.

The following theorem outlines a simulation of the first exit times of correlated drifted Brownian
motions.
\begin{theorem}
\label{thm:simulation3} Let $(\tau_1,\ldots,\tau_N)$ be the vector of first exit times
 given in Proposition \ref{prop:simulation} and denote by $f$ its
joint density. Let $(\chi_1^2,\ldots,\chi_N^2)$ be a vector of Chi-squared random variables with one degree of freedom satisfying
\begin{equation*}
\big(\chi_1^2,\ldots,\chi_N^2\big)\sim\big(H_1(\tau_1),\ldots,H_N(\tau_N)\big).
\end{equation*}
For $i=1,\ldots,N$,

\begin{itemize}
\item if $\mu_i\neq 0$, set
\begin{eqnarray}  \label{XX}
X_{i1}&=&\frac{b_i-x_i}{\mu_i}+\frac{\sigma_i^2\chi_i^2}{2\mu_i^2}-\frac{%
\sigma_i|\chi_i|}{2\mu_i^2} \sqrt{4\mu_i(b_i-x_i)+\sigma_i^2\chi_i^2};
\notag \\
X_{i2}&=&\frac{(b_i-x_i)^2}{\mu_i^2X_{i1}};
\end{eqnarray}

\item if $\mu_i=0$, set
\begin{equation}  \label{XX1}
X_{i1}=X_{i2}=\frac{(b_i-x_i)^2}{\sigma_i^2\chi_i^2}.
\end{equation}
\end{itemize}

For $k\in\{1,\ldots,N\}$ and $u_k\in\{1,2\}$, define $$S_u=\{1,2\}^N\setminus\{(u_1,\ldots,u_N)\}$$ and
\begin{eqnarray*}
&&p_{u_1u_2\cdots u_N}\\
&&=\bigg(1+%
\sum_{(i_1,\ldots,i_N)\in S_u}
\Big(\prod_{k=1}^N\Big|\frac{H_k'(X_{ku_k})} {H_k'(X_{ki_k})}\Big|\Big)\frac{%
f(X_{1i_1},X_{2i_2},\ldots,X_{Ni_N})}{f(X_{1u_1},X_{2u_2},\ldots,X_{Nu_N})}%
\bigg)^{-1},
\end{eqnarray*}
where $H_k'$ denotes the derivative of $H_k$. Denote by $\lambda$ the Lebesgue measure on $\mathbb R$.
Let $\{I_{u_1\cdots u_N}\}_{u_1,\ldots,u_N\in\{1,2\}}$ be any  partition of
interval $[0,1]$ satisfying
\begin{equation}
\label{I}\mathbb P\big(\lambda(I_{u_1\cdots
u_N})=p_{u_1\cdots u_N}\big)=1,
\end{equation} then we have
\begin{equation*}
(\tau_1,\ldots,\tau_N)\sim\sum_{u_1,\ldots,u_N\in\{1,2\}}\big(%
X_{1u_1},\ldots,X_{Nu_N}\big)\mathds{1}_{I_{u_1\cdots u_N}}(U),
\end{equation*}
where $U\sim Unif(0,1)$ is independent of $%
(\chi_1^2,\ldots,\chi_N^2)$.
\end{theorem}
 It is very difficult to generate the random vector $%
(\chi_1^2,\ldots,\chi_N^2)$ by its exact distribution since: the vector $(\chi_1^2,\ldots,\chi_N^2)$ can
not be exactly generated from a Gaussian copula. To show this fact we suppose $N=2$ and let $(Z_1,Z_2)$ be a Gaussian
vector verifying $Z_1^2=\chi_1^2$, $Z_2^2=\chi_2^2$, then we should have (see \cite{Peng}, Lemma $5.3.4$)
\begin{equation*}
Cov(\chi_1^2,\chi_2^2)=Cov(Z_1^2,Z_2^2)=2\big(Cov(Z_1,Z_2)\big)^2\geq 0.
\end{equation*}
But numerical integration (for instance by using MATLAB) shows $Cov(\chi_1^2,\chi_2^2)<0$ for some
values of $(x_1,b_1,\sigma_1,x_2,b_2,\sigma_2,\rho)$. For example, taking $%
b_1-x_1=b_2-x_2=-5,\sigma_1=\sigma_2=1,\rho=-0.5$ leads to $%
Cov(\chi_1^2,\chi_2^2)=-0.4007<0$.  Although it is difficult to actually generate two dependent first exit times when they have distributions that are not from a standard multivariate distribution, however, we can expect to simulate data following the same marginal
distribution and covariance structure via $(\tau_1,\ldots,\tau_N)$.  More precisely, we introduce the following approximation:

\begin{definition}
\label{defn:discov} Two second order (i.e., with finite covariance matrix) random
vectors $(X_1,\ldots,X_N)$ and $(Y_1,\ldots,Y_N)$ are said to be
approximately identically distributed if for $i\in\{1,\ldots,N\}$, their marginal distributions are equal: $X_i\sim Y_i$;
and their covariance matrices are identical:
\begin{equation*}
Cov\big((X_1,\ldots,X_N)\big)=Cov\big((Y_1,\ldots,Y_N)\big).
\end{equation*}
We denote this relation by $
(X_1,\ldots,X_N)\overset{approx}{\sim}(Y_1,\ldots,Y_N)$.
\end{definition}
Note that nowadays in financial risk aggregation,  elliptical copulas (including Gaussian copula) are widely used to simulate dependent data, since they provide a way to model correlated multivariate variables. In our case, we choose the Gaussian copula to generate first exit times following their exact means and covariance matrix. This is the so-called NORTA (normal-to-anything) method. Here we recall a nice result on the simulation of joint uniform random variables by following their covariance structure, starting from Gaussian vector:
\begin{lemma}
\label{lemme2} Assume that $U_1,U_2\sim Unif(0,1)$ have correlation $\rho_U$, then there exist $Z_1,Z_2\sim\mathcal N(0,1)$ with correlation $\rho_Z$
satisfying $
\rho_Z=2\sin\big(\frac{\pi}{6}\rho_U\big)$,
and
$$
(U_1,U_2)\overset{approx}{\sim} (\varphi(Z_1),\varphi(Z_2)),
$$
where $\varphi$ denotes the cumulative distribution function of the standard normal distribution.
\end{lemma}
This result is known as the exact relationship between Spearman
correlation and Bravais-Pearson correlation (see for example \cite{Hotelling}).  The approximation given in Definition \ref{defn:discov} and Lemma \ref{lemme2},  allow us  to simulate the first exit times in the following two situations:
\begin{description}
  \item[Case 1:] When $N=2$, and the Brownian motions have non-zero drifts.
  \end{description}
  In this case the transformations $\{H_k\}_{k=1,2}$ have multiple roots and the explicit formula of $f$, the joint density of $(\tau_1,\tau_2)$ can be derived. Hence one can determine the partitions $\{I_{ij}\}_{i,j=1,2}$ given in (\ref{I}) explicitly. The joint density $f$ of $(\tau_1,\tau_2)$ is given in the following:

      \begin{proposition}
\label{prop2} Let $(X_1(t),X_2(t))$ be a two-dimensional  Brownian
motion respectively with drifts $\mu_1,\mu_2\in\mathbb{R}$ and $
\sigma_1^2,\sigma_2^2>0$. Then starting from $(x_1,x_2)\in\mathbb{R}^2$, the
joint density of the first exit times $(\tau_1,\tau_2)$ to the barriers $
(b_1,b_2)\in\mathbb{R}^2$ is given as

\begin{enumerate}
\item Let $(\gamma_1,\gamma_2)=\big(\frac{\sigma_2\mu_1-\sigma_1\mu_2\rho}{
\sigma_1\sigma_2\sqrt{1-\rho^2}},\frac{\mu_2}{\sigma_2}\big)$, and $(\tilde{\mu_1
},\tilde{\mu_2})=\big((sgn(x_1-b_1))\gamma_1,(sgn(x_2-b_2))\gamma_2\big). $
For $0<s<t$,
\begin{eqnarray}  \label{density3}
&&f(s,t)=\sqrt{\frac{\pi}{2}}\frac{\sin\alpha}{\alpha^2s\sqrt{(t-s)^3}}e^{-r_0(\frac{r_0}{2s}+ \tilde{\mu_1}\cos\theta_0+\tilde{\mu_2}%
\sin\theta_0)-\frac{\tilde{\mu_1}^2s+\tilde{\mu_2}^2t}{2}}  \notag \\
&&\times\sum_{n=1}^{+\infty}n\sin\Big(\frac{n\pi(\alpha-\theta_0)}{\alpha}%
\Big)\int_0^{+\infty} e^{r\tilde{\mu_1}\cos\alpha-r^2(\frac{%
t-s\cos^2\alpha}{2s(t-s)})}I_{n\pi/\alpha}\Big(\frac{rr_0}{s}\Big)\, \mathrm{d}
r.\nonumber\\
\end{eqnarray}
\item For $0<t<s$,
\begin{eqnarray}  \label{density4}
&&f(s,t)=\sqrt{\frac{\pi}{2}}\frac{\sin\alpha}{\alpha^2t\sqrt{(s-t)^3}}\exp%
\Big(-r_0(\frac{r_0}{2t}+ \tilde{\mu_1}\cos\theta_0+\tilde{\mu_2}%
\sin\theta_0)  \notag \\
&&-\frac{(\tilde{\mu_1}^2+\tilde{\mu_2}^2)t}{2}-\frac{(\tilde{\mu}%
_1\sin\alpha -\tilde{\mu}_2\cos\alpha)^2(s-t)}{2}\Big)\sum_{n=1}^{+\infty}n%
\sin\Big(\frac{n\pi\theta_0}{\alpha}\Big)  \notag \\
&&\times\int_0^{+\infty} e^{-r(\tilde{\mu}_1\cos^2\alpha+\tilde{\mu%
}_2\sin\alpha\cos\alpha)-r^2\big(\frac{s-t\cos^2\alpha}{2t(s-t)}\big)%
}I_{n\pi/\alpha}\Big(\frac{rr_0}{t}\Big)\, \mathrm{d} r,  \notag \\
\end{eqnarray}
 where
       \begin{eqnarray*}
       \tilde{\rho}&=&\Big(sgn \big(\frac{b_1-x_1}{b_2-x_2}\big)\Big)\rho~\mbox{with $sgn(\cdot)$ being the sign function},\\
      \alpha &=& \left\{ \begin{array}{ll}
\pi+\tan^{-1}\big(-\frac{\sqrt{1-\tilde{\rho}^2}}{\tilde{\rho}}\big) & \textrm{if $\tilde{\rho}> 0$},\\
\frac{\pi}{2} &\textrm{if $\tilde{\rho}=0$},\\
\tan^{-1}\big(-\frac{\sqrt{1-\tilde{\rho}^2}}{\tilde{\rho}}\big) & \textrm{if $\tilde{\rho}<0$},
\end{array} \right.\\
        r_0&=&\frac{1}{\sigma_1\sigma_2}\sqrt{\frac{(b_1-x_1)^2\sigma_2^2+(b_2-x_2)^2\sigma_1^2-2|(b_1-x_1)(b_2-x_2)|\tilde{\rho}\sigma_1\sigma_2}{1-\tilde{\rho}^2}},\\
       \theta_0 &=& \left\{ \begin{array}{ll}
\pi+\tan^{-1}\big(\frac{\sigma_1|b_2-x_2|\sqrt{1-\tilde{\rho}^2}}{|b_1-x_1|\sigma_2-\tilde{\rho}|b_2-x_2|\sigma_1}\big) & \textrm{if $|b_1-x_1|\sigma_2<\tilde{\rho}|b_2-x_2|\sigma_1$},\\
\frac{\pi}{2} & \textrm{if $|b_1-x_1|\sigma_2=\tilde{\rho}|b_2-x_2|\sigma_1$},\\
\tan^{-1}\big(\frac{\sigma_1|b_2-x_2|\sqrt{1-\tilde{\rho}^2}}{|b_1-x_1|\sigma_2-\tilde{\rho}|b_2-x_2|\sigma_1}\big) & \textrm{if $|b_1-x_1|\sigma_2>\tilde{\rho}|b_2-x_2|\sigma_1$},\\
\end{array} \right.\\
\end{eqnarray*}
and $I_{\beta}$ with $\beta>0$ denotes the modified Bessel function of the
first kind of order $\beta$.
\end{enumerate}
\end{proposition}
 Inspired by Lemma \ref{lemme2}, we provide an approximation of the distribution of $(\tau_1,\tau_2)$ starting from a standard Gaussian vector $(Z_1,Z_2)$:

\begin{proposition}
\label{thm:simulation} Let $(\tau _{1},\tau _{2})$ be the first exit times
of $(X_{1}(t),X_{2}(t))$ with
correlation $\rho \in (-1,1)$. Set $(Z_{1},Z_{2})$ with
standard normal components satisfying:
\begin{eqnarray}
&&Corr(Z_{1},Z_{2})\nonumber\\
&&=2\sin \Big(\frac{\pi }{6}Corr\big(2\varphi \big(\frac{|\mu
_{1}\tau _{1}-(b_{1}-x_{1})|}{\sigma _{1}\sqrt{\tau _{1}}}\big)-1,2\varphi %
\big(\frac{|\mu _{2}\tau _{2}-(b_{2}-x_{2})|}{\sigma _{2}\sqrt{\tau _{2}}}%
\big)-1\big)\Big).\nonumber\\  \label{corr1}
\end{eqnarray}%
 Let
\begin{equation}
\label{chisquare}
(\chi _{1}^{2},\chi _{2}^{2})=\Big(\big(\varphi ^{-1}\big(\frac{\varphi
(Z_{1})+1}{2}\big)\big)^{2},\big(\varphi ^{-1}\big(\frac{\varphi (Z_{2})+1}{2}\big)\big)%
^{2}\Big).
\end{equation}%
Let $U\sim Unif(0,1)$ be independent of $\chi_1^2$, $\chi_2^2$. Then, the distribution of $(\tau _{1},\tau _{2})$ can be approximated by that of
\begin{equation}
\sum_{i,j\in \{1,2\}}\big(%
X_{1i},X_{2j}\big)\mathds{1}_{I_{ij}}(U),  \label{tau}
\end{equation}
in the sense that
\begin{equation*}
\label{hatchi}
\big(\chi_1^2,\chi_2^2\big)\overset{approx}{\sim }\big(H_1(\tau_1),H_2(\tau_2)\big),
\end{equation*}where for $i,j\in \{1,2\}$, $\big(X_{1i},X_{2j}\big)$ is defined in $(\ref{XX})$,  and $I_{ij}$ is given in $(\ref{I})$.
\end{proposition}
This proposition is a straightforward consequence of Lemma \ref{lemme2} and
Theorem \ref{thm:simulation3} by taking $N=2$. Hence we omit its proof. Note that in (\ref{corr1}) the correlation
$$
Corr\Big(2\varphi \Big(\frac{|\mu
_{1}\tau _{1}-(b_{1}-x_{1})|}{\sigma _{1}\sqrt{\tau _{1}}}\Big)-1,2\varphi %
\Big(\frac{|\mu _{2}\tau _{2}-(b_{2}-x_{2})|}{\sigma _{2}\sqrt{\tau _{2}}}%
\Big)-1\Big)
$$
can be computed numerically by using the joint density of $(\tau_1,\tau_2)$.  Although calculating the correlation by numerical integration in two dimensions is as difficult as calculating the mean exit time, we emphasize that this calculation can be done offline, while exit time simulation
usually is not. For instance, in Section \ref{Applications}, one is given a single
model of joint asset returns, for which the correlation $Corr(Z_1,Z_2)$ must be calculated once. However, several contingent claims can be written on these assets. The valuation of each of them requires an independent simulation, and
every simulation shares the same value of $Corr(Z_1,Z_2)$, whose calculation does not need to be repeated. More importantly, as we will show in the next subsection, when $N>2$ only pairwise correlations need to be numerically calculated.
\begin{description}
  \item[Case 2:] When $N\ge 2$, and the Brownian motions are without drifts.
  \end{description}In this case the joint density of $(\tau_1,\ldots,\tau_N)$ is generally unknown (or too complicated to get). We remark that it is not yet tractable in the literature to exactly simulate $(\tau_1,\ldots,\tau_N)$ when
its joint density is unknown. Fortunately, when all the drifts are vanishing, the transformations $H_i$'s have single roots, and as a consequence we have avoided the problem of selection among multiple roots of transform using the joint density of first exit times. An example of high dimensional simulation is given in Section \ref{high dimension}. We are
going to show the generating method for this case as well as the algorithms in the following sections.

\subsection{First Exit Times of Brownian Motions with Zero Drift}

\label{Nodrift} When $(\mu_1,\mu_2)=0$,
the reduced joint density of $f(s,t)$ in Proposition \ref{prop2} and the exact simulation of $(\tau_1,\tau_2)$ are given in the following Theorem \ref{thm1} (i). From Theorem \ref{thm1} (ii), we see $\{H_i\}_{i=1,2}$ reduce to a single root transformations.

\begin{theorem}
\label{thm1}

\begin{description}
\item[(i)] The random vector $(\tau_1,\tau_2)$ has the following joint
density:

\begin{enumerate}
\item For $0<s<t$,
\begin{eqnarray*}  \label{density1}
f(s,t)&=&\frac{\pi \sin\alpha}{2\alpha^2\sqrt{s(t-s\cos^2\alpha)}(t-s)} e^
{-\frac{r_0^2(t-s\cos2\alpha)}{2s((t-s)+(t-s\cos2\alpha))}}  \notag
\\
&\times&\sum_{n=1}^{\infty}n\sin\Big(\frac{n\pi(\alpha-\theta_0)}{\alpha}%
\Big)I_{n\pi/2\alpha} \Big(\frac{r_0^2(t-s)}{2s((t-s)+(t-s\cos2\alpha))}\Big)%
.  \notag \\
\end{eqnarray*}

\item For $0<t<s$,
\begin{eqnarray*}  \label{density2}
f(s,t)&=&\frac{\pi \sin\alpha}{2\alpha^2\sqrt{t(s-t\cos^2\alpha)}(s-t)} e^{
-\frac{r_0^2(s-t\cos2\alpha)}{2t((s-t)+(s-t\cos2\alpha))}}  \notag
\\
&\times&\sum_{n=1}^{\infty}n\sin\Big(\frac{n\pi\theta_0}{\alpha}\Big)%
I_{n\pi/2\alpha} \Big(\frac{r_0^2(s-t)}{2t((s-t)+(s-t\cos2\alpha))}\Big).
\end{eqnarray*}
\end{enumerate}

\item[(ii)] There exist  two Chi-squared random variables $\chi_1^2,\chi_2^2$
verifying
\begin{equation*}
\mathbb{E}\big(\chi_1^2\chi_2^2\big)=\frac{(b_1-x_1)^2(b_2-x_2)^2}{%
(\sigma_1\sigma_2)^2}\int_0^{+\infty}\int_0^{+\infty}\frac{f(s,t)}{st}\,
\mathrm{d} s\, \mathrm{d} t,
\end{equation*}
such that, with probability $1$,
 $
(\tau_1,\tau_2)=\Big(\frac{(b_1-x_1)^2}{\sigma_1^2 \chi_1^2},\frac{%
(b_2-x_2)^2}{\sigma_2^2 \chi_2^2}\Big).
$
\end{description}
\end{theorem}
Theorem \ref{thm1} is a straightforward consequence of Proposition \ref{prop2} and
Proposition \ref{prop:simulation}, where one takes $\mu_1=\mu_2=0$. Remark that
Metzler \cite{Metzler2} presented the joint density of $(\tau_1,\tau_2)$ only in the
case when $x_i>b_i=0$, $i=1,2$. Here we have generalized his result to any
real values of $(x_i,b_i)$. Moreover, this theorem could provide a simple
and fast simulation of the first exit times. We refer to the
following corollary:

\begin{corollary}
\label{cor1} There exists a Gaussian vector $(Z_1,Z_2)$ with $Z_1,Z_2\sim\mathcal N(0,1)$ and
\begin{equation*}  \label{corr3}
Corr(Z_1,Z_2)=2\sin\Big(\frac{\pi}{6}Corr\big(2\varphi\big(\frac{|b_1-x_1|}{%
\sigma_1\sqrt{\tau_1}}\big)-1, 2\varphi\big(\frac{|b_2-x_2|}{\sigma_2\sqrt{%
\tau_2}}\big)-1\big)\Big),
\end{equation*}
 such that the distribution of $(\tau_1,\tau_2)$ can be approximated by that of
\begin{equation*}
\bigg(\frac{(b_1-x_1)^2}{\sigma_1^2 (\varphi^{-1}(\frac{%
\varphi(Z_1)+1}{2}))^2}, \frac{(b_2-x_2)^2}{\sigma_2^2 (\varphi^{-1}(\frac{%
\varphi(Z_2)+1}{2}))^2}\bigg)
\end{equation*}
in the sense that
\begin{equation*}
\Big((\varphi^{-1}(\frac{%
\varphi(Z_1)+1}{2}))^2,(\varphi^{-1}(\frac{%
\varphi(Z_2)+1}{2}))^2\Big)\overset{approx}{\sim }\Big(\frac{\big(b_1-x_1\big)^2%
}{\sigma_1^2\tau_1},\frac{\big(b_2-x_2\big)^2}{%
\sigma_2^2\tau_2}\Big).
\end{equation*}
\end{corollary}
By using the property that the zero-mean Gaussian vector's distribution is
determined only by its covariance matrix, the result of Corollary \ref{cor1}
can be straightforwardly extended to high-dimensional correlated
Brownian motions with zero drift.

\begin{corollary}
\label{cor2} Let the integer $N\geq 2$, let $(Z_1,\ldots,Z_N)\sim%
\mathcal{N}(0,\Sigma)$ be any Gaussian vector satisfying,
\begin{equation*}
\Sigma=(r_{ij})_{N\times N},~\mbox{with}~r_{ii}=1 ~%
\mbox{for
$i\in\{1,2,\ldots,N\}$}
\end{equation*}
and for different $i,j\in\{1,2,\ldots,N\}$,
\begin{equation}  \label{rij}
r_{ij}=2\sin\Big(\frac{\pi}{6}Corr\big(2\varphi\big(\frac{|b_i-x_i|}{\sigma_i%
\sqrt{\tau_i}}\big)-1, 2\varphi\big(\frac{|b_j-x_j|}{\sigma_j\sqrt{\tau_j}}%
\big)-1\big)\Big),
\end{equation}
where each pair $(\tau_i,\tau_j)$ is with parameters $%
(x_i,b_i,\sigma_i,x_j,b_j,\sigma_j,\rho_{ij})$. Then the random vector $\big(\tau_1,\ldots,\tau_N\big)$ can be approximated by
\begin{equation}  \label{simuT1}
\bigg(\frac{%
(b_1-x_1)^2}{\sigma_1^2 (\varphi^{-1}(\frac{\varphi(Z_1)+1}{2}))^2},\ldots,
\frac{(b_N-x_N)^2}{\sigma_N^2 (\varphi^{-1}(\frac{\varphi(Z_N)+1}{2}))^2}%
\bigg)
\end{equation}
in the sense that
\begin{equation*}
\bigg(\big(\varphi^{-1}(\frac{%
\varphi(Z_i)+1}{2})\big)^2\bigg)_{i=1,\ldots,N}\overset{approx}{\sim }\bigg(\frac{\big(b_i-x_i\big)^2%
}{\sigma_i^2\tau_i}\bigg)_{i=1,\ldots,N}.
\end{equation*}
\end{corollary}
\section{Algorithms and Statistical Tests}
\subsection{Algorithms of First Exit Times}
\label{algorithms} Now we introduce a simulation algorithm of first exit
times for  Brownian motions. We are mainly interested in sampling the first exit times from correlated Brownian motions. The algorithms are given as follows:
\begin{description}
\item[Case 1:] Two-dimensional Brownian Motions with Drifts.
\end{description}
For a general two-dimensional vector of  Brownian motions, the
algorithm can be provided based on Proposition \ref{thm:simulation}.

\begin{enumerate}
\item Generate a Gaussian vector $(Z_1,Z_2)$ satisfying (\ref{corr1}).

\item For $i=1,2$, determine $(X_{i1},X_{i2})$ by using (\ref{XX}) and (\ref%
{XX1}).

\item Determine $(\chi_1^2,\chi_2^2)$ according to (\ref{chisquare}).

\item Generate $(\tau_1,\tau_2)$ by using (\ref{tau}).
\end{enumerate}
\begin{description}
\item[Case 2:] $N$-dimensional Brownian Motions without Drifts.
\end{description}
When the drift is zero, we supply an algorithm to simulate $%
(\tau_1,\ldots,\tau_N)$ based on Corollary \ref{cor2}.

\begin{enumerate}
\item Calculate $Corr\Big(2\varphi\big(\frac{|b_i-x_i|}{\sigma_i\sqrt{\tau_i}%
}\big)-1,2\varphi\big(\frac{|b_j-x_j|}{\sigma_j\sqrt{\tau_j}}\big)-1\Big)$
for $i,j\in\{1,\ldots,N\}$, $i\neq j$ by using $f_{ij}$, the joint density of $%
(\tau_i,\tau_j)$:
\begin{eqnarray*}
&&Corr\Big(2\varphi\big(\frac{|b_i-x_i|}{\sigma_i\sqrt{\tau_i}}\big)%
-1,2\varphi\big(\frac{|b_j-x_j|}{\sigma_j\sqrt{\tau_j}}\big)-1\Big) \\
&&=12\Big(\int_{\mathbb R_+^2}\big(2\varphi\big(\frac{|b_i-x_i|%
}{\sigma_i\sqrt{s}}\big)-1\big) \big(2\varphi\big(\frac{|b_j-x_j|}{\sigma_j%
\sqrt{t}}\big)-1\big)f_{ij}(s,t)\, \mathrm{d} s\, \mathrm{d} t-\frac{1}{4}\Big).
\end{eqnarray*}

\item Determine $r_{ij}$ for $i,j\in\{1,\ldots,N\}$, $i\neq j$, by using
Equation (\ref{rij}).

\item Generate $(Z_1,\ldots,Z_N)\sim\mathcal{N}(0,\Sigma)$, with $%
\Sigma=(r_{ij})_{N\times N}$.

\item Generate $(\tau_1,\ldots,\tau_N)$ by using (\ref{simuT1}).
\end{enumerate}

We emphasize that only pairwise correlations need to
be numerically calculated. Thus, one of the main advantages of our algorithm
is its little numerical complexity: rather than performing an $N$-dimensional
integration the algorithm performs a calculation of two-dimensional integrations.
\subsection{A Statistical Test of the Algorithms}
\label{K-S Test}
Although the simulation of the joint first exit times is approximate, it has been shown that the true distribution of the joint first times is quite similar to the one generated by using Gaussian copula. For example, Overbeck and Schmidt \cite{Overbeck} compared a Gaussian copula model to a calibration time changed model using the fair basket default swap spreads data  and showed that the results are quite close. From their numerical results, "the seed variance of fair the first-to-default spread is less than $0.10\%$ and much smaller for the other spreads"; McLeish \cite{McLeish} established an estimator of the maxima of two correlated Brownian motions using some normal vector, and he showed that the approximation is remarkably accurate. In addition, Metzler \cite{Metzler1} (Pages 50-51) has compared the true distribution of the first exit time of a 3-dimensional correlated Brownian motions (in this example the barriers are supposed to be equal) to the one from its corresponding Gaussian copula, and concludes that the first exit time in a "correlated Brownian drivers" model is quite similar to the one in a Gaussian copula model. Here we use another  statistical approach to compare our method to Monte Carlo's method. More comparisons using numerical results are made in the next section. Recall that the Kolmogorov-Smirnov test (K-S test) is used to check the equality of probability distributions. In order to compare random vectors we take a high dimensional K-S test (see for instance \cite{Lopes}), programmed in MATLAB. Let $(\tau_1,\tau_2)$ be the random vector generated by using Monte Carlo algorithm and $(\widetilde{\tau_1},\widetilde{\tau_2})$ be the one generated by using transformation with multiple roots.

We set the parameters for $i=1,2$, $(x_i,b_i,\sigma_i)=(\log(5),0,1)$ and the test's significant level $\alpha=0.01$. $10^5$ realizations are generated by using each method. Now we set up the hypothesis:
$$
\mathcal H_0:~(\tau_1,\tau_2)\sim(\widetilde{\tau_1},\widetilde{\tau_2})~\mbox{and}~\mathcal H_1:~(\tau_1,\tau_2)\nsim (\widetilde{\tau_1},\widetilde{\tau_2}).
$$
The results $\mathcal H_0$ presented in the following Table 1 show that we won't reject the fact that $(\tau_1,\tau_2)\sim(\widetilde{\tau_1},\widetilde{\tau_2})$ with confidence $0.99$ when $N=2$. However, the dispersion could be increasing when the dimension $N$ increases.
\begin{table}[H]
\begin{centering}
\begin{tabular}{|c|c|}
\hline
  Parameters  & K-S Test Result by MATLAB \tabularnewline
\hline
$\mu_1=\mu_2=0$,~$\rho=0.1$ & $\mathcal H_0$ \tabularnewline
\hline
$\mu_1=\mu_2=-0.05$,~$\rho=0.1$ & $\mathcal H_0$ \tabularnewline
\hline
$\mu_1=\mu_2=-0.05$,~$\rho=0.5$ & $\mathcal H_0$ \tabularnewline
\hline
\end{tabular}
\caption{K-S test with parameters $(x_i,b_i,\sigma_i)=(\log(5),0,1)$,~$\alpha=0.01$}
\par\end{centering}
\label{tab:Simulation-results-K-S-Test}
\end{table}
\section{Application to a Multiple-factor Model of Portfolio Default}

\label{Applications} The simulation of first exit times of Brownian motions has its interests in many credit risk applications. Let $X_{i}$ be the logarithm of the total value of the
assets of firm $i$, where $i=1,...,N$. Assume that $X_{i}$
satisfies the stochastic differential equation (\ref{dXk}). For $i\neq j$, we denote by $%
\rho _{ij}$ the correlation per unit of time between $X_{i}$ \ and $X_{j}$.

The substantial decrease of a firm's asset value is the main reason of the
default of the firm. Black and Cox \cite{Black} defined the threshold
value as the minimum asset value of the firm required by the debt covenants.
If $X_{i}$ falls to the threshold value $b_{i}$, the bond holder are
entitled to a 'deficiency claim' which can force the firm into bankruptcy. Zhou \cite{Zhou} defined the barrier $b_{i}$ as the logarithm of the sum of the short-term debt principal and one half of
long-term debt principal of a firm. $x_{i}$ is the logarithm of the current total asset values of
firms. With input parameters $(\mu_i,\sigma_i,x_i,b_i,\rho_{ij})$, we denote the output by $P_{i}$, the probability that $i$ firms out of $N$ in
a portfolio have defaulted by a time horizon $T$. For simplicity, the subsequent numerical simulations are based
on the assumption that the underlying Brownian motions have the same drift terms
and standard deviations. Since, for financial products, most correlations
lie between $-0.5$ and $0.5$, three different levels of correlations ( i.e., we denote by $\rho=\rho_{ij}=
-0.5,0.1,0.5$) are examined in the following tables. Zhou \cite{Zhou}'s parameters are
applied here: $(\sigma_i,x_i,b_i)=(1,\log(5),0)$. Over a short horizon, quick default events are rare
and the probabilities of multiple defaults could be extremely small, hence a long
horizon $T=10$ (years) is used to overcome this problem. We present the multiple default probabilities using two  different methods: our method and the Euler-based Monte Carlo algorithm. We compare them to the method in \cite{Zhou} (although our method is more general than Zhou's result, we compare our
method to his for accuracy), which is exact. In the following
tables, the step size of the Euler-based Monte Carlo method is equal to $0.0015625$ and the total
number of scenarios is $10^{6}$. Note that the Monte Carlo method is always an alternative, however it becomes very computationally intensive when a large number of scenarios need to be performed.

\subsection{Simulations of First Exit Times of Two-dimensional Correlated
Brownian Motions}

The examples without drift terms and with drift terms are both examined in
the following Tables $2$-$7$. Three different level correlations are presented to demonstrate how the correlation
of logarithm of asset value affects the default probability distribution. In
Table \ref{tab:Simulation-results-comparison-nodrift-1}, $\rho =0.1$ means
the movements of two firms' logarithm of asset value are weakly correlated with each
other. The fact that $\rho =0.5$ shows these movements
are strongly positively correlated. The simulation results illustrate the fact that the
stronger are the assets' correlations, the higher are the probabilities of multiple
defaults. Zhou \cite{Zhou} indicated that the default correlation and
the asset level correlation have the same sign which explains our results
here. For instance, a drop in one firm's asset value leads to a decrease in
another firm's asset value (which is closer to the default boundary), and
then leads to a rise of probability of both firms' defaults. If two firms'
asset values do not move towards the same direction, this behavior leads to a
drop in likelihood of multiple default events.

Compared to the analytical results in \cite{Zhou}, our
method is not completely unbiased, but the results are very promising when $\rho>0$. In fact from Tables 2-5, the estimated value of our method has an error rate less than $3.57\%$. When $\rho<0$, our numerical results (Tables 6-7) show that the error rate is less than $20\%$. However, in practice we rarely consider the case $\rho<0$, since it is quite artificial (see \cite{Overbeck}).


\begin{table}[H]
\begin{centering}
\begin{tabular}{|c|c|c|c|}
\hline
 & Our Method  & Zhou (2001) & Euler-based\tabularnewline
\hline
$P_{2}$ & $0.390521$ & $0.386337$ & $0.380781$\tabularnewline
\hline
$P_{1}$ & $0.440707$ & $0.448901$ & $0.451572$\tabularnewline
\hline
$P_{0}$ & $0.168782$ & $0.164761$ & $0.167657$\tabularnewline
\hline
\end{tabular}
\par\end{centering}
\caption{Two-dimensional simulation results with $\protect\rho=0.1,\protect\mu_{1}=\protect\mu_{2}=0$}
\label{tab:Simulation-results-comparison-nodrift-1}
\end{table}

\begin{table}[H]
\begin{centering}
\begin{tabular}{|c|c|c|c|}
\hline
 & Our Method  & Zhou (2001) & Euler-based\tabularnewline
\hline
$P_{2}$ & $0.445721$ & $0.446907$ & $0.440361$\tabularnewline
\hline
$P_{1}$ & $0.426332$ & $0.424764$ & $0.428239$\tabularnewline
\hline
$P_{0}$ & $0.127957$ & $0.128328$ & $0.131400$\tabularnewline
\hline
\end{tabular}
\par\end{centering}
\caption{Two-dimensional simulation results with $\protect\rho=0.1,\protect\mu_{1}=\protect\mu_{2}=-0.05$}
\label{tab:Simulation-results-drift-1}
\end{table}

\begin{table}[H]
\begin{centering}
\begin{tabular}{|c|c|c|c|}
\hline
& Our Method  & Zhou (2001) & Euler-based\tabularnewline
\hline
$P_{2}$ & $0.439642$ & $0.445308$ & $0.440013$\tabularnewline
\hline
$P_{1}$ & $0.344621$ & $0.330958$ & $0.331241$\tabularnewline
\hline
$P_{0}$ & $0.215747$ & $0.223732$ & $0.228756$\tabularnewline
\hline
\end{tabular}
\par\end{centering}
\caption{Two-dimensional simulation results with $\protect\rho=0.5,\protect\mu_{1}=\protect\mu_{2}=0$}
\label{tab:Simulation-results-comparison-nodrift-2}
\end{table}
\begin{table}[H]
\begin{centering}
\begin{tabular}{|c|c|c|c|}
\hline
 & Our Method  & Zhou (2001) & Euler-based\tabularnewline
\hline
$P_{2}$ & $0.505348$ & $0.502006$ & $0.498137$\tabularnewline
\hline
$P_{1}$ & $0.30397$ & $0.314566$ & $0.315721$\tabularnewline
\hline
$P_{0}$ & $0.190682$ & $0.183426$ & $0.186142$\tabularnewline
\hline
\end{tabular}
\par\end{centering}
\caption{Two-dimensional simulation results with $\protect\rho=0.5,\protect\mu_{1}=\protect\mu_{2}=-0.05$}
\label{tab:Simulation-results-drift-2}
\end{table}
\begin{table}[H]
\begin{centering}
\begin{tabular}{|c|c|c|c|}
\hline
& Our Method  & Zhou (2001) & Euler-based\tabularnewline
\hline
$P_{2}$ & $0.325874$ & $0.308726$ & $0.301252$\tabularnewline
\hline
$P_{1}$ & $0.570430$ & $0.604123$ & $0.607972$\tabularnewline
\hline
$P_{0}$ & $0.103696$ & $0.087150$ & $0.090786$\tabularnewline
\hline
\end{tabular}
\par\end{centering}
\caption{Two-dimensional simulation results with $\protect\rho=-0.5,\protect\mu_{1}=\protect\mu_{2}=0$}
\label{tab:Simulation-results-comparison-nodrift-3}
\end{table}
\begin{table}[H]
\begin{centering}
\begin{tabular}{|c|c|c|c|}
\hline
 & Our Method  & Zhou (2001) & Euler-based\tabularnewline
\hline
$P_{2}$ & $0.372292$ & $0.376896$ & $0.367961$\tabularnewline
\hline
$P_{1}$ & $0.566252$ & $0.564787$ & $0.571982$\tabularnewline
\hline
$P_{0}$ & $0.061456$ & $0.058316$ & $0.060057$\tabularnewline
\hline
\end{tabular}
\par\end{centering}
\caption{Two-dimensional simulation results with $\protect\rho=-0.5,\protect\mu_{1}=\protect\mu_{2}=-0.05$}
\label{tab:Simulation-results-drift-3}
\end{table}

\subsection{First Exit Times Simulations for High-dimensional Correlated
Brownian Motions}
\label{high dimension}
We note
that our method is roughly as good as Euler-based Monte Carlo method. This is
important, as the method in \cite{Zhou} is not applicable in higher dimensions. Some simulation results of first exit times for high-dimensional Brownian motions  are listed in Table \ref%
{tab:Simulation-results-comparison-three-dim}. The subsequent simulations are
based on the simple examples without drift terms. Since Zhou's results can not
work on higher dimensional cases, only Monte Carlo method is used here. For simplicity, all the correlations are set to be $0.1$. For instance,
the correlation matrix $Corr(X_{1},X_{2},X_{3})$ of three firms' logarithm of asset value $X_{i}$ where $i=1,2,3$ in
Table \ref{tab:Simulation-results-comparison-three-dim} is
\begin{eqnarray*}
\begin{array}{rcl}
Corr(X_{1},X_{2},X_{3}) & = & \left(
\begin{array}{rcl}
1 & 0.1 & 0.1 \\
0.1 & 1 & 0.1 \\
0.1 & 0.1 & 1%
\end{array}%
\right).
\end{array}%
\end{eqnarray*}

The simulation results are consistent with the fact that the Monte Carlo
method tends to under-estimate the multiple defaults probabilities.
According to \cite{Carmona},
when the number of firms increases, this estimation could have larger deviation. However, compared to the results in \cite{Zhou} which can not
calculate multiple default probabilities for more than two firms, the new
proposed method can be advantageously expanded to high dimensional Brownian motions. In Table 8, the error rate of our estimation compared to Monte Carlo algorithm is less than $5.38\%$, which is acceptable in practice.

\begin{center}
\begin{table}[H]
\begin{centering}
\begin{tabular}{|c|c|c|}
\hline
 & Our Method & Euler-based\tabularnewline
\hline
$P_{3}$ & $0.257971$ & $0.250633$\tabularnewline
\hline
$P_{2}$ & $0.395472$ & $0.403186$\tabularnewline
\hline
$P_{1}$ & $0.267637$ & $0.271008$\tabularnewline
\hline
$P_{0}$ & $0.078920$ & $0.075173$\tabularnewline
\hline
\end{tabular}
\par\end{centering}
\caption{Three-dimensional simulation results with $Corr(X_{1},X_{2},X_{3})$}
\label{tab:Simulation-results-comparison-three-dim}
\end{table}
\end{center}
\vspace{-2cm}
\section{Conclusion}

A new numerical algorithm is proposed to solve the first exit times of a vector of  Brownian motions with zero and non-zero drifts. Compared to the methods in \cite{Metzler1,Metzler2} and Zhou's analytical and numerical results, the advantages of our method are the following: When $N=2$, our algorithms can take into account drifts and can be easily extended to high dimensions, especially when the drifts are zero. As a consequence, our method allows calculating the expected value of more general functions of the first exit times, not constrained to the minimum of the latter,
which is the traditional advantage of simulation.

\section{Appendix}

\subsection{Proof of Theorem \protect\ref{thm:simulation3}.}
 Let $(\tau_1,\tau_2,\ldots,\tau_N)$ be a vector of first
exit times with  parameters $\big\{(\mu_i,\sigma_i,x_i,b_i)\big\}_{i=1,\ldots,N}$. We assume that $x_i\neq b_i$ and $\frac{\mu_i}{%
b_i-x_i}\ge 0$ for $i\in\{1,\ldots,N\}$ and the joint density of $%
(\tau_1,\ldots,\tau_N)$ is $f$. Hence, still by \cite{Shuster}, $(\tau_1,\ldots,\tau_N)$ is the solution of the following
equation:
\begin{equation}
\label{chi}
\big(H_1(\tau_1),\ldots,H_N(\tau_N)\big)=\big(\chi_1^2,\ldots,\chi_N^2\big).
\end{equation}
Notice that for $i=1,\ldots,N$, $v_i\ge0$, the equation
$$
\big(H_1(x_1),\ldots,H_N(x_n)\big)=(v_1,\ldots,v_N)
$$
has $2^N$ roots (can be surplus) $\big\{\big(x_{1u_1},x_{2u_2},%
\ldots,x_{Nu_N}\big):~(u_1,\ldots,u_N)\in\{1,2\}^N\big\}$,

\begin{itemize}
\item if $\mu_i\neq 0$, set
\begin{eqnarray*}
x_{i1}&=&\frac{b_i-x_i}{\mu_i}+\frac{\sigma_i^2v_i}{2\mu_i^2}-\frac{\sigma_i%
}{2\mu_i^2} \sqrt{4\mu_i(b_i-x_i)v_i+\sigma_i^2v_i^2};  \notag \\
x_{i2}&=&\frac{(b_i-x_i)^2}{\mu_i^2x_{i1}};
\end{eqnarray*}

\item if $\mu_i=0$, set
\begin{equation}  \label{X1}
x_{i1}=x_{i2}=\frac{(b_i-x_i)^2}{\sigma_i^2v_i}.
\end{equation}
\end{itemize}

The remaining problem is choosing one among the $2^N$ roots according to
the observed value of the Chi-squared vector $\big(\chi_1^2,\ldots,\chi_N^2\big)$. We explain how this
process can be taken. For a particular root $\big(x_{1u_1},x_{2u_2},%
\ldots,x_{Nu_N}\big)$, we compute the probability of choosing it. Let $%
\epsilon>0$ be arbitrarily small and let $H_i^{-1}(v_i-\epsilon,v_i+%
\epsilon) $ be the inverse image of $(v_i-\epsilon,v_i+\epsilon)$. For $%
i=1,\ldots,N$, let the intervals $(y_{iu_i}^{(1)},y_{iu_i}^{(2)})\subset
H_i^{-1}(v_i-\epsilon,v_i+\epsilon)$ be disjoint for each $u_i\in\{1,2\}$
and contain the $u_i$-th root of equation $H_i(x)=v_i$: $x_{iu_i}%
\in(y_{iu_i}^{(1)},y_{iu_i}^{(2)})$. Notice that we should choose $\big(%
x_{1u_1},x_{2u_2},\ldots,x_{Nu_N}\big)$ as the observation of $%
(\tau_1,\ldots,\tau_N)$ if $\tau_i\in(y_{iu_i}^{(1)},y_{iu_i}^{(2)})$ for $%
i=1,\ldots,N$. Therefore the probability of choosing the particular root $%
\big(x_{1u_1},\ldots,x_{Nu_N}\big)$ is given as
\begin{eqnarray*}
&&\mathbb{P}^{\epsilon}\big(\big(x_{1u_1},\ldots,x_{Nu_N}\big)\big)\\
&&=\frac{\mathbb{P}\big(\tau_1\in(y_{1u_1}^{(1)},y_{1u_1}^{(2)}),\ldots,%
\tau_N\in(y_{Nu_N}^{(1)},y_{Nu_N}^{(2)})\big)}{\sum_{j_1,\ldots,j_N\in\{1,2%
\}}\mathbb{P}\big(\tau_i\in(y_{1j_1}^{(1)},y_{1j_1}^{(2)}),\ldots,\tau_N%
\in(y_{Nj_N}^{(1)},y_{Nj_N}^{(2)})\big)}.
\end{eqnarray*}
Observe that $(y_{iu_i}^{(1)},y_{iu_i}^{(2)})\rightarrow x_{iu_i}$ as $%
\epsilon\rightarrow0^+$. Hence letting $\mathbb{P}^{\epsilon}$ denote the
probability measure which depends on $\epsilon$ and letting $%
\epsilon\rightarrow0^+$, by a similar principle to L'H\^opital's rule (see Equation (3) in \cite{Michael}), we get
\begin{eqnarray*}
&&\mathbb{P}_{(u_1,u_2,\ldots,u_N)}((v_1,\ldots,v_N)):=\lim_{\epsilon%
\rightarrow0^+}\mathbb{P}^{\epsilon}\big(\big(x_{1u_1},\ldots,x_{Nu_N}\big)%
\big) \\
&&=\lim_{\epsilon\rightarrow0^+}\bigg\{1+\sum_{(j_1,\ldots,j_N)\in S_u} \frac{\mathbb P\big(\tau_1\in(y_{1j_1}^{(1)},y_{1j_1}^{(2)}),\ldots,%
\tau_N\in(y_{Nj_N}^{(1)},y_{Nj_N}^{(2)})\big)}{\mathbb P\big(\tau_1%
\in(y_{1u_1}^{(1)},y_{1u_1}^{(2)}),\ldots,\tau_N%
\in(y_{Nu_N}^{(1)},y_{Nu_N}^{(2)})\big)}\bigg\}^{-1} \\
&&=\bigg\{1+\sum_{(j_1,\ldots,j_N)\in S_u}\lim_{\epsilon\rightarrow0^+} \frac{\frac{\mathbb P%
\big(\tau_1\in(y_{1j_1}^{(1)},y_{1j_1}^{(2)}),\ldots,\tau_N%
\in(y_{Nj_N}^{(1)},y_{Nj_N}^{(2)})\big)}{(y_{1j_1}^{(2)}-y_{1j_1}^{(1)})%
\cdots(y_{Nj_N}^{(2)}-y_{Nj_N}^{(1)})}}{\frac{\mathbb P\big(\tau_1%
\in(y_{1u_1}^{(1)},y_{1u_1}^{(2)}),\ldots,\tau_N%
\in(y_{Nu_N}^{(1)},y_{Nu_N}^{(2)})\big)}{(y_{1u_1}^{(2)}-y_{Nu_N}^{(1)})%
\cdots(y_{Nu_N}^{(2)}-y_{Nu_N}^{(1)})}}\\
&&~~\times \prod_{i=1}^N\bigg(\frac{\frac{%
(y_{ij_i}^{(2)}-y_{ij_i}^{(1)})}{%
2\epsilon}}{\frac{(y_{iu_i}^{(2)}-y_{iu_i}^{(1)})}{2\epsilon}}\bigg)\bigg\}^{-1} \\
&&=\bigg\{1+\sum_{(j_1,\ldots,j_N)\in S_u} \frac{f\big(%
x_{1j_1},\ldots,x_{Nj_N}\big)}{f\big(x_{1u_1},\ldots,x_{Nu_N}\big)}\prod_{i=1}^N\Big|%
\frac{H_i^{\prime }(x_{iu_i})}{H_i^{\prime
}(x_{ij_i})}\Big|\bigg\}^{-1}, \\
\end{eqnarray*}
where we recall that $S_u=\{1,2\}^N\setminus\{(u_1,\ldots,u_N)\}$ and for $i=1,\ldots,N$ and any $x>0$,
 $$
H_i^{\prime }(x)=\frac{((\mu_ix)^2-(b_i-x_i)^2)}{\sigma_i^2x^2}.
$$
Let $(v_1,\ldots,v_N)$ be equal to the
observation of Chi-squared vector defined in (\ref{chi}) and let $%
(X_{1u_1},\ldots,X_{Nu_N})$ be the root labeled $(u_1,\ldots,u_N)$ of (\ref{chi}), then
the 'probability' (it is a random variable) of choosing this root is given as $$
p_{u_1\ldots u_N}=\mathbb{P}_{(u_1,\ldots,u_N)}((\chi_1^2,\ldots,\chi_N^2)).$$
It remains to generate the probability distribution $\{p_{u_1\ldots
u_N}\}_{(u_1,\ldots,u_N)\in\{1,2\}^N}$. It can be simulated by generating an
independent uniform random variable $U\sim Unif(0,1)$ and the following random partition of interval $[0,1]$:
 $$
[0,1]=\bigcup_{u_1,\ldots,u_N\in\{1,2\}}I_{u_1\ldots u_N},
$$
with for almost every $\omega\in\Omega$, $\mathbb{P}(U\in I_{u_1\ldots u_N}(\omega))=p_{u_1\ldots u_N}(\omega)$. This is equivalent to (\ref{I}). Finally $$\big(H_1(X_{1u_1}),\ldots,H_N(X_{Nu_N})\big)=(\chi_1^2,\ldots,\chi_N^2)$$
and
\begin{equation*}
(\tau_1,\ldots,\tau_N)\sim
\sum_{u_1,\ldots,u_N\in\{1,2\}}(X_{1u_1},\ldots,X_{Nu_N})\mathds{1}_{
I_{u_1\ldots u_N}}(U).~\square
\end{equation*}

\subsection{Proof of Proposition \protect\ref{prop2}.}

Notice that we derive the idea from the seminal work by Iyengar \cite{Iyengar}. Unfortunately, this work contains errors. Metzler \cite{Metzler2} provided the correct formula of the two-dimensional first exit times joint density, however the one with non-zero drift with drift is not explicitly given. In this section we provide a relatively closed form of the joint density. Our work also extends the joint
density of first exit times in \cite{Metzler2} to the case where the barriers could be any real
values. First we assume that for $i=1,2$, $x_i-b_i>0$. The main idea is to
transform the  Brownian motions to independence (see \cite{Iyengar}). By such a linear transformation, the first exit time remains
invariant when the Brownian motion starts from $(x-b)$ to the barrier $0$. Hence,
without loss of generality, we only consider the barriers as horizontal
axis $Y=0$. Define $T:\mathbb{R}^2\rightarrow \mathbb{R}^2$ by
\begin{equation*}
T(x)=\Big({}_0^{\sigma_1\sqrt{1-\rho^2}}~{}_{\sigma_2}^{\sigma_1\rho}\Big)%
^{-1}x,
\end{equation*}
 be the transformation of the vector of
Brownian motions $X(t)=(X_1(t),X_2(t))$ with correlation $\rho$ and variances $\sigma_1^2,\sigma_2^2$ to the independent standard Brownian motions. Denote the latter by $%
Z(t)=(Z_1(t),Z_2(t))=T(X(t))$ and $z_0=Z(0)=(r_0\cos
\theta_0,r_0\sin\theta_0)$. After transformation $T$, the horizontal barrier line $Y=0$ remains
the same and $Y=b_2$ turns to be the line $Y=(\tan\alpha) X$. Note that the
notations $(r_0,\theta_0,\alpha)$ are exactly the same as mentioned in the
joint density of $(\tau_1,\tau_2)$. The explicit form of joint density depends on the conditions $%
\tau_1<\tau_2$ and $\tau_2<\tau_1$.

\begin{enumerate}
\item For $0< s<t$, this shows $\tau_1<\tau_2$. Following the argument of Metzler \cite{Metzler2}, the joint density of $(\tau_1,\tau_2)$ is given as
\begin{eqnarray*}
&&\int_{r=0}^{+\infty}\frac{\mathbb{P}(\tau_1\in\, \mathrm{d} s,
Z(\tau_1)\in\, \mathrm{d} z,z=r(\cos\alpha,\sin\alpha))}{\, \mathrm{d} s}
\\
&&~~\times\frac{\mathbb{P}(\tau_2-\tau_1\in\, \mathrm{d} (t-s)|\tau_1\in\,
\mathrm{d} s, Z(\tau_1)\in\, \mathrm{d} z,z=r(\cos\alpha,\sin\alpha))}{%
\, \mathrm{d} (t-s)} \\
&&=\int_0^{+\infty}\!\!\!\!\!e^{\gamma^{\prime }(z-z_0)-\frac{|\gamma|^2s}{2}%
}\frac{\pi}{\alpha^2sr}e^{-\frac{r^2+r_0^2}{2s}}\sum_{n=1}^{+%
\infty} n\sin(\frac{n\pi(\alpha-\theta_0)}{\alpha})I_{n\pi/\alpha}(\frac{rr_0%
}{s}) \\
&&~~\times\Big(\frac{r\sin\alpha}{\sqrt{2\pi(t-s)^3}}e^{{-\frac{%
(r\sin\alpha+\gamma_2(t-s))^2}{2(t-s)}}}\Big)\, \mathrm{d} r,
\end{eqnarray*}
where
\begin{description}
\item[(i)]
$
\gamma=(\gamma_1~\gamma_2)^{\prime }=T\Big({}_{\mu_2}^{\mu_1}\Big)=\Big(%
\frac{\sigma_2\mu_1-\sigma_1\mu_2\rho}{\sigma_1\sigma_2\sqrt{1-\rho^2}}~%
\frac{\mu_2}{\sigma_2}\Big)^{\prime }$,
with $\gamma^{\prime }$ being the transpose of $\gamma$ and $|\gamma|$ being
the Euclidian norm.
\item[(ii)] By easy computation, we can get $\alpha$, $r_0$, $\theta_0$ as in Proposition \ref{prop2}, with $\tilde\rho=\rho$.
\item[(iii)] $z=(r\cos\alpha,r\sin\alpha)$ is the polar coordinates of the
transformation of the exit position $Z(\tau_1)=(Z_1(\tau_1),Z_2(\tau_1))$.

\item[(iv)] By Iyengar \cite{Iyengar} and Metzler \cite{Metzler2},
\begin{eqnarray*}
&&\frac{\mathbb{P}(\tau_1\in\, \mathrm{d} s, Z(\tau_1)\in\, \mathrm{d}
z,z=r(\cos\alpha,\sin\alpha))}{\, \mathrm{d} s} \\
&&=e^{\gamma^{\prime }(z-z_0)-\frac{|\gamma|^2s}{2}}\frac{\pi%
}{\alpha^2sr}e^{-\frac{r^2+r_0^2}{2s}}\sum_{n=1}^{+\infty} n\sin(\frac{%
n\pi(\alpha-\theta_0)}{\alpha})I_{n\pi/\alpha}(\frac{rr_0}{s})
\end{eqnarray*}
is the density of time for $Z(t)$ to first exit from the line $Y=(\tan\alpha) X$ and
\begin{eqnarray*}
&&\frac{\mathbb{P}(\tau_2-\tau_1\in\, \mathrm{d} (t-s)|\tau_1\in\,
\mathrm{d} s, Z(\tau_1)\in\, \mathrm{d} z,z=r(\cos\alpha,\sin\alpha))}{%
\, \mathrm{d} (t-s)} \\
&&=\frac{r\sin\alpha}{\sqrt{2\pi(t-s)^3}}e^{{-\frac{%
(r\sin\alpha+\gamma_2(t-s))^2}{2(t-s)}}}\, \mathrm{d} r
\end{eqnarray*}
is in fact the inverse Gaussian density on $t-s$ which denotes the remaining
time for $Z_2(t)$ to first exit its barrier, once $Z_1(t)$ exits. By martingale
property, this event can be regarded as the first exit time of a standard
 Brownian motion with drift $\gamma_2$, starting from $r\sin\alpha$%
, to the barrier $0$ (see Fig. $1$ for $\tau_1<\tau_2$).
\end{description}

\item For $0< t<s$ (this shows $\tau_2<\tau_1$), the joint density can be
obtained in a similar way:
\begin{eqnarray*}
&&\frac{\mathbb{P}(\tau_2\in\, \mathrm{d} t, Z(\tau_2)\in\, \mathrm{d}
z,z=(r,0))}{\, \mathrm{d} t}\\
&& =\int_{r=0}^{+\infty}\frac{\mathbb{P}(\tau_2\in\, \mathrm{d} t,
Z(\tau_2)\in\, \mathrm{d} z,z=(r,0))}{\, \mathrm{d} t}\\
&&~~\times\frac{\mathbb{P}(\tau_1-\tau_2\in\, \mathrm{d} (s-t)|\tau_2\in\,
\mathrm{d} t, Z(\tau_2)\in\, \mathrm{d} z,z=(r,0))}{\, \mathrm{d} (s-t)}\\
&&=\int_0^{+\infty}e^{\gamma^{\prime }(z-z_0)-\frac{|\gamma|^2t}{2}}%
\frac{\pi}{\alpha^2tr}e^{-\frac{r^2+r_0^2}{2t}}\sum_{n=1}^{+\infty}
n\sin(\frac{n\pi\theta_0}{\alpha})I_{n\pi/\alpha}(\frac{rr_0}{t})\\
&&~~\times\Big(\frac{r\sin\alpha}{\sqrt{2\pi(s-t)^3}}e^{-\frac{%
(r\sin\alpha+(s-t)(\gamma_1\cos(\alpha-\pi/2)+\gamma_2\cos(\pi-\alpha)))^2}{%
2(s-t)}}\Big)\, \mathrm{d} r,
\end{eqnarray*}
where we just need to remark that $z=(r,0)$ and the latter inverse Gaussian
density (on $s-t$) denotes the remaining time for $Z_1(t)$ to first exit its
barrier, once $Z_2(t)$ exits. It can be further regarded as the first exit
time of a standard  Brownian motion with drift $%
\gamma_1\cos(\alpha-\pi/2)+\gamma_2\cos(\pi-\alpha))$, starting from $%
r\sin\alpha$, until the barrier $0$.
\end{enumerate}

Now it suffices to simplify the above two formulae. For example, let us
simplify the first one: when $0<s<t$,
\begin{eqnarray*}
&&f(s,t)=\int_0^{+\infty}e^{\gamma_1(r\cos\alpha-r_0\cos\theta_0)+%
\gamma_2(r\sin\alpha-r_0\sin\theta_0)-\frac{(\gamma_1^2+\gamma_2^2)s}{2}}
\\
&&\times\frac{\pi}{\alpha^2sr}e^{-\frac{r^2+r_0^2}{2s}}%
\sum_{n=1}^{+\infty}n\sin\big(\frac{n\pi(\alpha-\theta_0)}{\alpha}\big)%
I_{n\pi/\alpha}(\frac{rr_0}{s}) \\
&&\times\frac{r\sin\alpha}{\sqrt{2\pi(t-s)^3}}e^{-\frac{%
(r\sin\alpha+\gamma_2(t-s))^2}{2(t-s)}}\, \mathrm{d} r \\
&=&\frac{\pi\sin\alpha}{\alpha^2s\sqrt{2\pi(t-s)^3}}e^{%
-\gamma_1r_0\cos\theta_0-\gamma_2r_0\sin\theta_0-\frac{r_0^2}{2s}-\frac{%
\gamma_1^2s+\gamma_2^2t}{2}} \\
&&\times\int_0^{+\infty}\sum_{n=1}^{+\infty}n\sin\big(\frac{%
n\pi(\alpha-\theta_0)}{\alpha}\big)e^{r\gamma_1\cos\alpha-r^2\big(%
\frac{1}{2s}+\frac{s\sin^2\alpha}{2s(t-s)}\big)}I_{n\pi/\alpha}(\frac{%
rr_0}{s})\, \mathrm{d} r \\
&=&\sqrt{\frac{\pi}{2}}\frac{\sin\alpha}{\alpha^2s\sqrt{(t-s)^3}}e^{%
-r_0(\frac{r_0}{2s}+ \gamma_1\cos\theta_0+\gamma_2\sin\theta_0)-\frac{%
\gamma_1^2s+\gamma_2^2t}{2}} \\
&&\times\sum_{n=1}^{+\infty}n\sin\big(\frac{n\pi(\alpha-\theta_0)}{\alpha}%
\big)\int_0^{+\infty} e^{r\gamma_1\cos\alpha-r^2(\frac{%
t-s\cos^2\alpha}{2s(t-s)})}I_{n\pi/\alpha}(\frac{rr_0}{s})\, \mathrm{d}
r.
\end{eqnarray*}
When $0<t<s$, we can pursue the same approach to simplify the expression
(see Fig. $1$ for $\tau_1>\tau_2$). Finally, we get

\begin{enumerate}
\item For $0<s<t$,
\begin{eqnarray}  \label{jointdensity1}
&&f(s,t)=\sqrt{\frac{\pi}{2}}\frac{\sin\alpha}{\alpha^2s\sqrt{(t-s)^3}}e^
{-r_0(\frac{r_0}{2s}+ \gamma_1\cos\theta_0+\gamma_2\sin\theta_0)-\frac{%
\gamma_1^2s+\gamma_2^2t}{2}}  \notag \\
&&~~\times\sum_{n=1}^{+\infty}n\sin\big(\frac{n\pi(\alpha-\theta_0)}{\alpha}%
\big)\int_0^{+\infty} e^{\gamma_1\cos\alpha-r^2(\frac{%
t-s\cos^2\alpha}{2s(t-s)})}I_{n\pi/\alpha}(\frac{rr_0}{s})\, \mathrm{d}
r.\nonumber\\
\end{eqnarray}
\item For $0<t<s$,
\begin{eqnarray}  \label{jointdensity2}
&&f(s,t)=\sqrt{\frac{\pi}{2}}\frac{\sin\alpha}{\alpha^2t\sqrt{(s-t)^3}}\exp%
\Big(-r_0\big(\frac{r_0}{2t}+ \gamma_1\cos\theta_0+\gamma_2\sin\theta_0\big)  \notag
\\
&&~~-\frac{(\gamma_1^2+\gamma_2^2)t}{2}-\frac{(\gamma_1\sin\alpha
-\gamma_2\cos\alpha)^2(s-t)}{2}\Big)\sum_{n=1}^{+\infty}n\sin\big(\frac{%
n\pi\theta_0}{\alpha}\big)  \notag \\
&&~~\times\int_0^{+\infty} e^{-r\big(\gamma_1\cos^2\alpha+\gamma_2\sin\alpha%
\cos\alpha\big)-r^2\big(\frac{s-t\cos^2\alpha}{2t(s-t)}\big)}I_{n\pi/\alpha}(%
\frac{rr_0}{t})\, \mathrm{d} r.
\end{eqnarray}
\end{enumerate}
\begin{figure}[H]
\vspace{-4cm}
~~~~~~~~~~\includegraphics[width=11cm]{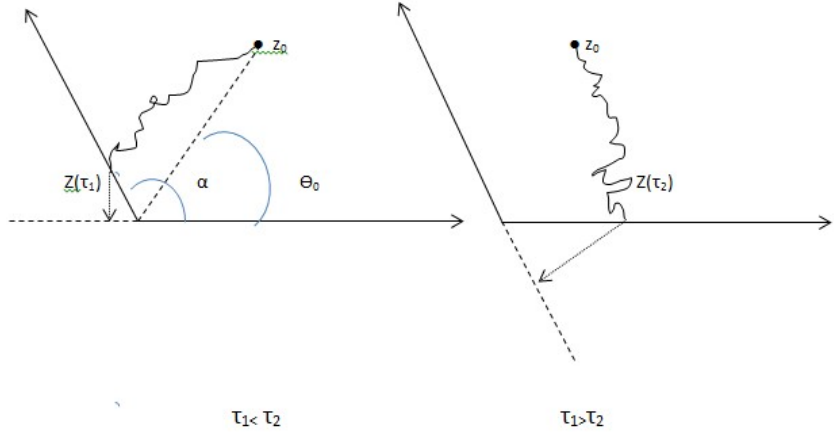}
\vspace{-4cm}
\caption{An example of the trajectories of $\{Z(t)\}_{t\ge0}$ for $\protect\tau_1<\protect\tau_2$ and $\protect\tau_1>\protect\tau_2$}
\end{figure}
Now we are going to derive the most generalized density, where the barriers
could be any real values. To this end we first observe that, when $x_1-b_1>0$, $x_2-b_2>0$, there exists a function $\tilde f$ such that
\begin{equation*}
f(s,t)=\tilde{f}\big(x_1-b_1,x_2-b_2,\mu_1,\mu_2,\sigma_1,\sigma_2,\rho,s,t\big)
\end{equation*}
By using the symmetric property of Brownian motion, when $x_i-b_i<0$ for $%
i=1 $ or $2$, we notice that
\begin{eqnarray*}
\tau_i&=&\inf\big\{t> 0: X_i(t)=b_i,~X_i(0)=x_i\big\}\\
&=&\inf\big\{t>0:-X_i(t)=-b_i,-X_i(0)=-x_i\big\},
\end{eqnarray*}
where $\{-X_i(t)\}_{t\geq 0}$ is a Brownian motion with drift $%
-\mu_i=(sgn(x_i-b_i))\mu_i$ starting from $-x_i$. Therefore
the first exit times of $(X_1(t),X_2(t))$ with drifts $(\mu_1,\mu_2)$, starting
from $(x_1,x_2)$ to the barriers $(b_1,b_2)$ are almost surely equal to the first exit
times of $\big((sgn(x_1-b_1))X_1(t),(sgn(x_2-b_2))X_2(t)\big)$, starting
from the initial values $\big((sgn(x_1-b_1))x_1,(sgn(x_2-b_2))x_2\big)$ to the barriers $%
((sgn(x_1-b_1))b_1,(sgn(x_2-b_2))b_2)$. Remark that
\begin{equation*}
Corr\big((sgn(x_1-b_1))X_1(t),(sgn(x_2-b_2))X_2(t)\big) =\frac{sgn(b_1-x_1)}{%
sgn(b_2-x_2)}\rho=\tilde{\rho}.
\end{equation*}
Finally, we have for all $x_i\neq b_i$, $i=1,2$,
\begin{equation*}
f(s,t)=\tilde{f}\big(|b_1-x_1|,|b_2-x_2|,(sgn(x_1-b_1))\mu_1,(sgn(x_2-b_2))\mu_2,\sigma_1,\sigma_2,
\tilde{\rho},s,t\big).
\end{equation*}
For any $x_i\neq b_i$, $i=1,2$, plugging the arguments $$%
(|b_1-x_1|,|b_2-x_2|,(sgn(x_1-b_1))\mu_1,(sgn(x_2-b_2))\mu_2,\sigma_1,\sigma_2,\tilde{\rho})$$
into (\ref{jointdensity1}) and (\ref{jointdensity2}), we get Proposition \ref%
{prop2}. $\square$

\subsection{Proof of Corollary \protect\ref{cor1}.}

To prove Corollary \ref{cor1}, we rely on Lemma \ref{lemme2} and the following lemma which is known as a consequence of Sklar's Theorem:

\begin{lemma}
\label{lemme:covZ:} Let $F_{\chi^2}$ be the cumulative probability
distribution function of $\chi^2\sim \chi^2(1)$, then $U=F_{\chi^2}(\chi^2)\sim Unif(0,1)$.
\end{lemma}
Now we are ready to prove Corollary \ref{cor1}. Since
\begin{equation*}
(\chi_1^2,\chi_2^2)=\Big(\frac{(b_1-x_1)^2}{\sigma_1^2\tau_1},\frac{%
(b_2-x_2)^2}{\sigma_2^2\tau_2}\Big),
\end{equation*}
then by Lemma \ref{lemme:covZ:},
\begin{eqnarray*}  \label{cov1}
(U_1,U_2):=\Big(F_{\chi^2}\big(\frac{(b_1-x_1)^2}{\sigma_1^2\tau_1}\big)%
,F_{\chi^2}\big(\frac{(b_2-x_2)^2}{\sigma_2^2\tau_2}\big)\Big)
\end{eqnarray*}
has marginal distribution $Unif(0,1)$. Observe that for any $x\geq 0$,
$
F_{\chi^2}(x)=\mathbb{P}(Z^2\leq x)=2\varphi(\sqrt{x})-1
$.
Then the correlation of $(U_1,U_2)$ can be given by using the joint density
of $(\tau_1,\tau_2)$:
\begin{eqnarray*}
&&Corr(U_1,U_2)=Corr\big(F_{\chi^2}\big(\frac{(b_1-x_1)^2}{\sigma_1^2}\big),F_{\chi^2}%
\big(\frac{(b_2-x_2)^2}{\sigma_2^2}\big)\big) \\
&&=Corr\big(2\varphi(\frac{|b_i-x_i|}{\sigma_i\sqrt{\tau_i}})-1,2\varphi(%
\frac{|b_j-x_j|}{\sigma_j\sqrt{\tau_j}})-1\big) \\
&&=\frac{\mathbb{E}\big(\big(2\varphi(\frac{|b_i-x_i|}{\sigma_i\sqrt{\tau_i}}%
)-1\big)\big(2\varphi(\frac{|b_j-x_j|}{\sigma_j\sqrt{\tau_j}})-1\big)\big)-%
\mathbb{E}(U_1)\mathbb{E}(U_2) }{\sqrt{Var(U_1)Var(U_2)}} \\
&&=12\Big(\int_0^{+\infty}\!\!\!\!\int_0^{+\infty}\!\!\!\!\big(2\varphi(\frac{|b_i-x_i|}{%
\sigma_i\sqrt{s}})-1\big) \big(2\varphi(\frac{|b_j-x_j|}{\sigma_j\sqrt{t}})-1%
\big)f_{ij}(s,t)\, \mathrm{d} s\, \mathrm{d} t-\frac{1}{4}\Big).
\end{eqnarray*}
Therefore by Lemma \ref{lemme2}, there exists a Gaussian vector such that (%
\ref{corr3}) holds. By the fact that for $i=1,2$,
\begin{equation*}
\chi_i^2=F_{\chi^2}^{-1}(U_i)\sim F_{\chi^2}^{-1}(\varphi(Z_i)),
\end{equation*}
and by observing that $F_{\chi^2}^{-1}(x)=\big(%
\varphi^{-1}(\frac{x+1}{2})\big)^2$, Corollary \ref{cor1} holds. $\square$




\begin{thebibliography}{99}

\bibitem{Anonymous} Anonymous (2009). In defense of the Gaussian copula. {\it Economist, April 29.}

\bibitem{Black} Black, F. and Cox, J. C. (1976). Valuing corporate securities: some effects
of bond indenture provisions. {\it Journal of Finance, {\bf 31} (2), 351-367.}

\bibitem{Broadie} Broadie, M., Glasserman, P. and Kou, S. (1997).  A continuity correction
for discrete barrier options. {\it Mathematical Finance, {\bf 7} (4), 325-348.}

\bibitem{Carmona} Carmona, R., Fouque, J. P. and Vestal, D. (2009). Interacting particle
systems for the computation of rare credit portfolio losses. {\it Finance
and Stochastics, {\bf 13} (4), 613-633.}

\bibitem{Collin-Dufresne} Collin-Dufresne, P. and Goldstein, R. S. (2001). Do credit spreads reflect stationary leverage ratios? {\it The Journal of Finance, {\bf 56} (5), 1929-1957.}

\bibitem{Fahim} Fahim, A., Touzi, N. and Warin, X. (2011). A probabilistic numerical
method for fully nonlinear parabolic PDEs. {\it The Annals of Applied
Probability, {\bf 21} (4), 1322-1364.}

\bibitem{Glasserman} Glasserman, P. and Li, J. (2005). Importance sampling for portfolio credit
risk. {\it Management Science, {\bf 51} (11), 1643-1656.}

\bibitem{Hotelling} Hotelling, H. and Pabst, M. R. (1936). Rank correlation and tests of significance
involving no assumption of normality. {\it Annals of Mathematical
Statistics, {\bf 7} (1), 29-43.}

\bibitem{Huh} Huh, J. and Kolkiewicz, A. (2008). Computation of multivariate barrier
crossing probability and its applications in credit risk models. {\it North
American Actuarial Journal, {\bf 12} (3), 263-291.}

\bibitem{Iyengar} Iyengar, S. (1985). Hitting lines with two-dimensional Brownian motion.
{\it SIAM Journal on Applied Mathematics, {\bf 45} (6), 983-989.}

\bibitem{Lopes} Lopes, R. H. C. (2007). The two-dimensional Kolmogorov-Smirnov test. {\it XI International Workshop on Advanced Computing and Analysis Techniques in Physics Research, Proceedings of Sciences.}

\bibitem{Metzler1} Metzler, A. (2008). Multivariate first passage models in credit risk. {\it Ph. D.
Dissertation, University of Waterloo.} https://uwspace.uwaterloo.ca/bitstream/handle/10012/\\4090/thesis$\_$adam$\_$metzler.pdf?sequence=1

\bibitem{Metzler2}  Metzler, A. (2010). On the first passage problem for correlated Brownian motion.
{\it Statistics and Probability Letters, {\bf 80} (5), 277-284.}

\bibitem{McLeish} McLeish, D. L. (2004). Estimating the correlation of processes using extreme values. {\it Fields Institute Communications, {\bf 44}, 447-467.}

\bibitem{Michael} Michael, J. R., Schucany, W. R. and Haas, R. W. (1976). Generating random
variates using transformations with multiple roots, {\it The American
Statistician, {\bf 30} (2), 88-90.}

\bibitem{Milstein} Milstein, G. N. and Tretyakov, M. V. (1999). Simulation of a space-time
bounded diffusion, {\it Annals of Applied Probability, {\bf 9} (3), 732-779.}

\bibitem{Muller} Muller, M. E. (1956). Some continuous Monte Carlo methods for the Dirichlet
problem. {\it The Annals of Mathematical Statistics, {\bf 27} (3), 569-589.}

\bibitem{Overbeck} Overbeck, L. and Schmidt, W. (2005). Modeling default dependence with threshold models. {\it Journal of Derivatives, Summer 10-19.}

\bibitem{Peng} Peng, Q. (2011). Statistical inference for hidden multifractional processes
in a setting of stochastic volatility models. {\it Ph. D. Dissertation, University
Lille $1$.} https://ori-nuxeo.univ-lille1.fr/nuxeo/site/esupversions/c56bf4f8-263c-4dc2-9196-85e5ef5d22e1


\bibitem{Shevchenko} Shevchenko, P. V. (2003). Addressing the bias in Monte Carlo pricing of
multi-asset options with multiple barriers through discrete sampling.
{\it The Journal of Computational Finance, {\bf 6}, 1-20.}

\bibitem{Shuster} Shuster, J. (1968). On the inverse Gaussian distribution function. {\it Journal
of the American Statistical Association, {\bf 63} (324), 1514-1516.}

\bibitem{Zhou} Zhou, C. (2001). An analysis of default correlations and multiple defaults,
{\it The Review of Financial Studies, {\bf 14} (2), 555-576.}
\end{thebibliography}
\end{document}